\documentclass[12pt]{article}

\def\W{{\mbox{\mcal W}\kern1pt}}

 \font\msyma = msam10
\def\leq{\mathrel{\mbox{\msyma 6}}}
\def\geq{\mathrel{\mbox{\msyma >}}}

\font\msymb = msbm10
\def\leftroundarrow{{\hbox{\msymb\char"78}}}
\def\rightroundarrow{{\hbox{\msymb\char"79}}}

\newcommand\mexa{\nopagebreak \begin{flushleft}\smallskip \nopagebreak
                \begin{minipage}[c]{1cm}\sloppy}

\newcommand\mexc{\end{minipage}\kern -3cm \smallskip\end{flushleft}}

\newcounter{subs}

\def\unt0x{\int \limits_{0}^{\scriptstyle x}}
\def\i0z{\int \limits_{0}^{\scriptstyle \zeta}}

\def\lprodint#1#2{\int\limits_{#1}^{{{\kern5pt{\scriptstyle #2}}\atop\leftroundarrow}}}
\def\rprodint#1#2{\int\limits_{#1}^{{{\kern4pt{\scriptstyle #2}}\atop\rightroundarrow}}}

\begin{document}
\begin{center}\textbf{Integrable Operators and Canonical
Differential Systems}\end{center}

\begin{center}\textbf{Lev Sakhnovich}\end{center}

My address: Lev Sakhnovich, 735, Crawford, Brooklyn, 11223,New
York, NY, USA.\\
My tel. (718) 646-2757.\\
My e-mail address:Lev.Sakhnovich@verizon.net\\
Research Professor,University of Connecticut at Storrs.
\newpage

\begin{abstract}
In this article we consider a class of integrable operators and
investigate its connections with the following theories:the
spectral theory of non-self-adjoint operators, the Riemann-Hilbert
problem, the canonical differential systems and the random
matrices theory.
\end{abstract}

\section*{Introduction}
In the article [12] we considered the operators of the type
\begin{equation}
Sf=L(x)f(x)+P.V.\int_{a}^{b}\frac{D(x,t)}{x-t}f(t)dt,\end{equation}
where $f(x){\in}L_{k}^{2}(a,b)$ and $k{\times}k$ matrix functions
$L(x)$ and $D(x,t)$ are such that \begin{equation}
L(x)=L^{\star}(x),\quad D(x,t)=-D^{\star}(t,x).\end{equation} (The
symbol P.V. indicates that the corresponding integral is
understood as the principal value.) \\
Later in the work [8] the important class of the operators S ,when
\begin{equation} k=1, \quad L(x)=1,\quad D(x,x)=0 ,
\end{equation} was studied in details. These results had a number
of interesting applications [5],[8].\\
In our works [12],[13]  the connection of the operators S with the
spectral theory of non-selfadjoint operators was shown.The
operator identity \begin{equation}
(QS-SQ)f=\int_{a}^{b}D(x,t)f(t)dt, \quad Qf=xf(x),\end{equation}
plays an essential role in these articles. From the identity (4)
follows the statement.\\
 \textbf{Proposition 1.} \emph{Let the
kernel $D(x,t)$be degenerate , i.e. $D(x,t)=iA(x)A^{\star}(t)$,
where $A(x)$ is a $k{\times}m$ matrix function ($k{\leq}m$). If
the operator $S$ is invertible , then the operator $T=S^{-1}$ has
the form \begin{equation}
Tf=M(x)f(x)+P.V.\int_{a}^{b}\frac{E(x,t)}{x-t}f(t)dt,\end{equation}
where $ M(x)=M^{\star}(x)$  and the kernel $E(x,t)$ is also
degenerate and has the form \begin{equation}
E(x,t)=iB(x)B^{\star}(t) ,\end{equation} $B(x)$ is a $k{\times}m$
matrix function.}\\
The operators $S$ and $T$ lead to the Riemann-Hilbert matrix
problem \begin{equation}
W_{+}(\sigma)=W_{-}(\sigma)R^{2}(\sigma),\quad
a{\leq}\sigma{\leq}b,\end{equation} where $m{\times}m$ matrix
function $W(z)$ is analytic , when $z{\notin}[a,b]$. Here matrix
function $R^{2}(\sigma)$ can be constructed with the help of the
operators $S$ and $T,\quad W_{\pm}(\sigma)$ is defined by the
relation
\begin{equation}
 W_{\pm}(\sigma)=\lim W(z), \quad y{\to}0,\quad z={\sigma}+iy.
 \end{equation}
 In the present article an essential role is played by the
 canonical differential system \begin{equation}
 \frac{d}{dx}W(x,z)=i\frac{JH(x)}{z-x}W(x,z), \quad
 W(0,z)=I_{m},\end{equation} where $m{\times}m$ matrix $J$ is such
 that\\
 $J=J^{\star},\quad J^{2}=I_{m}$ and $H(x){\geq}0$.\\
 The monodromy matrix of system (9) coincides with the solution of
 the Riemann-Hilbert problem (7),  i.e.
 \begin{equation}
 W(z)=W(b,z) .\end{equation}
 Let us note that $W(z)$ is a characteristic matrix function of
 the operator (see [2],[10])
 \begin{equation}
 Af=xf+i\int_{a}^{x}\beta(x)J\beta^{\star}(t)f(t)dt,\quad f(x){\in}L_{k}^{2}(a,b), \end{equation}
 where $\beta(x)$ is a $k{\times}m$ matrix function such, that
 \begin{equation}
 \beta^{\star}(x)\beta(x)=H(x).\end{equation} We deduce in this
 article a new sufficient condition of the linear similarity of
 the operator $A$ to the operator $Qf=xf$. It easily follows from
 (9)  that $W(x,z)$ in the neighborhood of $z=\infty$ admits the
 representation \begin{equation}
 W(x,z)=I_{m}+\frac{M_{1}(x)}{z}+\frac{M_{2}(x)}{z^{2}}+...,
 \end{equation} where \begin{equation}
 M_{1}(x)=i\int_{a}^{x}JH(t)dt.\end{equation}
 In view of (9) and (14) all the coefficients $M_{k}(x)$ are defined
 if the coefficient $M_{1}(x)$ is known. This fact is of interest
 as the representation \begin{equation}
 W(b,z)=I_{m}+\frac{M_{1}(b)}{z}+\frac{M_{2}(b)}{z^{2}}+...
 \end{equation} is closely connected with the problems of the
 random matrices theory [4],[14]. From the view point of the
 random matrix theory it is important that in this article the
 procedure of constructing the matrix $M_{1}(x)$ is given (section
 3). We pay the principal attention to the matrix
 version of the class (3), when \begin{equation}
 k{\geq}1,\quad L(x)=I_{k},\quad D(x,x)=0.\end{equation}
 For this class the corresponding matrix function $R^{2}(x)$ from
 (7) has a special structure, namely  \begin{equation}
 [R^{2}(x)-I_{m}]^{2}=0.\end{equation}
 The corresponding matrix function $JH(x)$ is nilpotent when
 $m=1$, i.e.
 \begin{equation}
 [JH(x)]^{2}=0.\end{equation} In the last part of the paper we
 consider a number of examples.
\section{Integrable operators and Riemann-Hilbert problem}
In this section we remind of a number of facts contained in the
paper [12]. We use these facts in the next sections.
Let $W(z)$ be $m{\times}m$ matrix function.\\
We suppose that the following conditions are fulfilled.\\
1). Matrix function $W(z)$ is analytic in the domain
$z{\notin}[a,b],(-\infty<a<b<\infty)$ and satisfies the equality
\begin{equation}
 W(z)=I_{m}+\frac{1}{2{\pi}i}\int_{a}^{b}\frac{F(x)}{x-z}dx,
\end{equation}
where $F(x)$ is bounded $m{\times}m$ matrix function on the
segment [a,b].\\
 2). The relations
\begin{equation} W^{\star}(z)JW(\bar{z})=J ,
\end{equation}
\begin{equation}
 i\frac{W^{\star}(z)JW(z)-J}{z-\bar{z}}{\geq}0 ,\quad
 z{\ne}\bar{z}
\end{equation}
are true. ( Here $m{\times}m$ matrix J satisfies the equalities
 $J=J^{\star},\quad J^{2}=I$).\\
The equality (1) guarantees the almost everywhere existence of the
limits
\begin{equation}
W_{\pm}(x)=limW(z)\quad as \quad y{\to}{\pm}0 ,\quad z=x+iy.
\end{equation}
Now we use the polar decomposition (see [11])
\begin{equation}
W_{+}(x)=U(x)R(x) ,\quad W_{-}(x)=U(x)R^{-1}(x) ,
\end{equation}
where $m{\times}m$ matrix functions $U(x)$ and $R(x)$ are such
that
\begin{equation}
U^{\star}(x)JU(x)=J ,\quad JR(x)=R^{\star}(x)J
\end{equation}
and in addition the spectrum of $R(x)$ is positive.\\
Matrix function $R(x)$ is called $J$-\emph{module} of matrix
function $W_{+}(x)$. By relations (23) and (24) we have
\begin{equation}
R^{2}(x)=JW_{+}^{\star}(x)JW_{+}(x).
\end{equation}
According to the theory of J-module [11] the relations
\begin{equation}
D(x)=J[R(x)-R^{-1}(x)]{\geq}0 ,\quad x{\in}[a,b] ,
\end{equation}
\begin{equation}
D(x)=0 , \quad x{\notin}[a,b]
\end{equation}
are true. Now we introduce the matrix functions $F_{1}(x),
F_{2}(x)$   with the help of the relations
\begin{equation}
F_{1}^{\star}(x)F_{1}(x)=D(x),\quad
F_{2}(x)=F_{1}(x)JU^{\star}(x).
\end{equation}
\emph{Remark} 1. Matrix functions $F_{1}(x)$ and $F_{2}(x)$ are
$k{\times}m$ matrices, where $k=\sup[\ rank D(x)],
a{\leq}x{\leq}b$.
Hence $k{\leq}m$.\\
Using relations (23),(26) and (28) we can write
 \begin{equation}
W_{+}(x)-W_{-}(x)=F_{2}^{\star}(x)F_{1}(x)=F(x) .\end{equation} In
addition to conditions 1) and 2) we suppose:\\
3). The matrix functions $ F_{1}(x)$ and
$F_{2}(x)$are bounded on segment [a,b].\\
 Let us define the
operators $\Pi$ and $\Gamma$ by
formulas ${\Pi}g=\frac {1}{\sqrt{2\pi}}F_{1}(x)g$ ,\\
${\Gamma}g=-\frac {i}{\sqrt{2\pi}}F_{2}(x)g,$ where $g$ are
$m{\times}1$ vectors, ${\Pi}g$ and ${\Gamma}g$ belong to
$L_{k}^{2}(a,b)$. Then we have
\begin{equation}
{\Pi}^{\star}f(x)=\frac{1}{\sqrt{2\pi}}{\int}_{a}^{b}F_{1}^{\star}(x)f(x)dx,
\end{equation}
\begin{equation}
{\Gamma}^{\star}f(x)=\frac
{i}{\sqrt{2\pi}}{\int}_{a}^{b}F_{2}^{\star}(x)f(x)dx,\end{equation}
where $f(x){\in}L_{k}^{2}(a,b)$. The next assertion follows from
formulas (19) ,(30) and (31).\\
\textbf{Proposition 2}. \emph{The matrix function W(z) admits the
realization \begin{equation}
W(z)=I_{m}-{\Gamma}^{\star}(Q-zI)^{-1}{\Pi},\end{equation} where
the operator $Q$ is defined by the relation}
\begin{equation} Qf=xf,\quad
f(x){\in}L_{k}^{2}(a,b).\end{equation} Next we introduce the
$k{\times}k$ matrix \begin{equation}
L(x)=[I_{k}+\frac{1}{4}(F_{1}(x)JF_{1}^{\star}(x))^{2}]^{1/2}
\end{equation} and consider the operators
\begin{equation}
Sf=L(x)f(x)+\frac{i}{2\pi}P.V.\int_{a}^{b}\frac{F_{1}(x)JF_{1}^{\star}(t)}{x-t}f(t)dt
,
\end{equation}
\begin{equation}
Tf=L(x)f(x)-\frac{i}{2\pi}P.V.\int_{a}^{b}\frac{F_{2}(x)JF_{2}^{\star}(t)}{x-t}f(t)dt.
\end{equation}
 The introduced operators S and T are acting
in the space $L_{k}^{2}(a,b)$ and $f(x)$ is a $k{\times}1$ vector
function.\\
\textbf{Theorem 1}.(see [13], p.45-46) \emph{The operators S and T
are positive , bounded and}
\begin{equation}
T=S^{-1},\quad SF_{2}(x)=F_{1}(x)J.
\end{equation}
From relation (23) we deduce that \begin{equation}
W_{+}(x)=W_{-}(x)R^{2}(x),\quad x{\in}[a,b]\end{equation}
\begin{equation}
W_{+}(x)=W_{-}(x),\quad x{\notin}[a,b]\end{equation}
 Formulas
(38) and (39) lead to the
Riemann-Hilbert Problem.\\
 \textbf{Problem }1. \emph{To recover the
matrix function $W(z)$ by the given J-module}
$R(x)$.\\
In the case $J=I$ Problem 1 plays an essential role in the
prediction theory of the stationary processes [15]. The case when
$J{\ne}I$ is
important for the theory of random matrices [5], [8],[14].\\
We  solve Problem 1 in the following way.\\
1.By the given matrix $R^{2}(x)$ we construct the matrix $D(x)$
(see (26)).\\
2.Using the first of equalities (28) we find $F_{1}(x)$.\\
3.With the help of formula (1) the operator $S$ is constructed.\\
4.Due to  the second equality of (37) we have
$F_{2}(x)=S^{-1}F_{1}(x)J$.\\
5.Now it is easy to see that formulas (19) and (29) give the
solution of the Riemann-Hilbert problem (7) with the normalizing
condition
\begin{equation}
W(z){\to}I \quad as \quad z{\to}\infty. \end{equation}
 \emph{Remark 2}. The
operators $S$ and $T$ defined by formulas (1) and (5) are called
integrable [5], [8].The case when $k=1$ and
\begin{equation}
F_{1}(x)JF_{1}^{\star}(x)=0 \end{equation} has  important
applications in the theory of the random matrices(see [4], [7],
[8], [14]). The general case was used in the spectral theory of
the non-selfadjoint operators [12],[13].
\section{Spectral theory}
We introduce some important notions .\\
Let the linear bounded operator have the form
\begin{equation}
A=A_{R}+iA_{I}, \end{equation} where $A_{R}$ and $A_{I}$ are
self-adjoint operators acting in Hilbert space $H$.There is a
bounded linear operator $K$ which maps a Hilbert space $G$ in $H$
so that \begin{equation} A_{I}=KJK^{\star},\end{equation} where
$J$ acts in $G$ and $J=J^{\star},\quad J^{2}=I$.\\
\textbf{Definition} 2 (see [2], [10)].\emph{The operator function
\begin{equation} W(\lambda)=I-2iK^{\star}(A-{\lambda}I)^{-1}KJ
\end{equation} is
called the characteristic operator function of $A$.}\\
We recall that the simple part of $A$ means the operator  which is
induced by $A$ on the subspace
$H_{1}=\overline{\sum_{k=0}^{\infty}A^{k}D_{A}}$, where
$D_{A}=\overline{(A-A^{\star})H}$. In paper [12] we deduced
Theorem 1 for the case $m{\leq}\infty$. From this fact we obtain
the following assertion [12],[13].\\
\textbf{Theorem} 2. \emph{If the characteristic operator function
$W(z)$ of the operator A satisfies the condition }\begin{equation}
||W(z)||{\leq}c,\quad z{\ne}\bar{z} \end{equation}\emph{ for some
c , then the simple part of A is linearly similar to a
self-adjoint
operator with a absolutely continuous spectrum}\\
It follows from relation (45) that $W(z)$ satisfies the conditions
1)-3).The converse is not true. Using this fact we receive a new
version of Theorem 2.\\
\textbf{Theorem} 3.  \emph{If the characteristic operator function
$W(z)$ of the operator $A$ satisfies the conditions 1)-3) , then
the statement of Theorem 2 is true.}\\
\textbf{Example}. We consider the case when  \begin{equation}
F_{1}(x)=[x+i,x-i],\quad 0{\leq}x{\leq}1,\quad
j=\left[\begin{array}{cc}
  -1 & 0 \\
  0 & 1
\end{array}\right]\end{equation}. The corresponding operator $S$
has the form
\begin{equation}
Sf=f(x)-\frac{1}{\pi}\int_{0}^{1}f(t)dt.\end{equation} Due to
relations (46) and (47) we have \begin{equation}
F_{2}(x)=[-q(x),\overline{q(x)}], \end{equation} where
\begin{equation}
q(x)=x+\frac{1}{2(\pi-1)}+i\frac{\pi}{\pi-1}. \end{equation} Using
the property of the  Cauchy integral  (see[6]) we deduce from
relation (19) that \begin{equation}
W(z)=-\frac{1}{2{\pi}i}F(0)logz +0(1),\quad z{\ne}\bar{z},\quad
|z|<\frac{1}{2},\end{equation}
\begin{equation} W(z)=-\frac{1}{2{\pi}i}F(1)log(z-1)
+0(1),\quad z{\ne}\bar{z},\quad |z-1|<\frac{1}{2}.\end{equation}
It follows from formulas (46) and (48),(49) that $F(0){\ne}0,\quad
F(1){\ne}0.$ Hence the constructed $W(z)$ satisfies the conditions
of Theorem 3 but does not satisfy the condition (45) of Theorem 2.
\section{Canonical differential systems}
It follows from Theorem 3 that the  following operator
\begin{equation}
S_{\xi}f=L(x)f(x)+\frac{i}{2\pi}P.V.\int_{a}^{\xi}\frac{F_{1}(x)JF_{1}^*(t)}{x-t}f(t)dt
\end{equation}
 is positive , bounded and invertible.\\
 We set
 \begin{equation}
 \Phi(\xi,x)=S_{\xi}^{-1}F_{1}(x),
\end{equation}
\begin{equation}
B(\xi)=\frac{1}{2\pi}\int_{a}^{\xi} \Phi^*(\xi,x)F_{1}(x)dx .
\end{equation}
\textbf{Lemma} 1. \emph{The matrix function $B(\xi)$ is absolutely
continuous and monotonically
increasing .}\\
\emph{Proof}. As it is known [3],[9] the operator $S^{-1}$ can be
represented in the form \begin{equation} S^{-1}= U^{\star}U ,
\end{equation}where the linear bounded operator $U$  acts in the
space $L_{k}^{2}(a,b)$ and satisfies the condition
\begin{equation}
U^{\star}P_{\xi}=P_{\xi}U^{\star}P_{\xi},\quad
a{\leq}{\xi}{\leq}b,\end{equation} where $P_{\xi}f(x)=f(x),
a{\leq}x{\leq}{\xi}$  and  $P_{\xi}f(x)=0,\quad {\xi}{\leq}x.$
From relations (54) and (55) we deduce the equality
\begin{equation}
\frac{d}{dx}
B(x)=\emph{H}(x)=\frac{1}{2\pi}h^{\star}(x)h(x),\end{equation}
where
\begin {equation}
h(x)=UF_{1}(x).\end{equation} The lemma is proved.\\
 Let us consider the system of equations
\begin{equation}
W(x,z)=I+iJ\int_{a}^{x}\frac{dB(\xi)}{z-\xi}W(\xi,z).
\end{equation}
\textbf{Theorem 4.} (see[13], Ch.3) \emph{The following equality
\begin{equation}W(b,z)=W(z)
\end{equation} holds.}

\textbf{Corollary 1}. \emph{The integral system (59) is equivalent
to the differential system }\begin{equation}
\frac{dW(x,z)}{dx}=\frac{iJ\emph{H}(x)}{z-x}W(x,z),\quad
H(x){\geq}0
\end{equation} \emph{with the boundary condition $W(a,z)=I_{m}$.
Here the matrix function
H(x) is defined by relation (57) .} \\
\textbf{Corollary 2}. \emph{The matrix function $W(z)$ is the
monodromy
matrix of system (61), i.e. $W(z)=W(b,z)$.}\\
 Due to (61)
in the neighborhood of $z=\infty$ the following relation
\begin{equation} W(x,z)= I+M_{1}(x)/z+M_{2}(x)/z^{2}+...\end{equation}
is fulfilled. It follows from (59) and (61) that
\begin{equation}
 M_{1}(x)=iJB(x).
\end{equation}
Formulas (53) (54) and (63) give the solution of the following
inverse problem.\\
\textbf{Problem 2}. \emph{To recover the matrix function $H(x)$
and $M_{1}(x)$ by the given J-module $R(x)$.}
Theorem 3 and relation (54)  imply the following assertion.\\
\textbf{Proposition 3}. \emph{ If equality \begin{equation}
F_{1}(x)=0, \quad \alpha {\leq}x{\leq}\beta,\quad
\alpha{\ne}\beta\end{equation} is true  then}
\begin{equation}
F_{2}(x)=0,\quad W_{+}(x)=W_{-}(x), \quad R(x)=I,\quad
\alpha{\leq}x{\leq}\beta.\end{equation} \textbf{Corollary 3}.
 \emph{If condition (64) is fulfilled then} \begin{equation}
B^{\prime}(x)=\emph{H}(x)=0,\quad \alpha
{\leq}x{\leq}\beta.\end{equation}
\section{Examples}
\textbf{Example 1}. Let us consider the case when \begin{equation}
J=j=\left[\begin{array}{cc}
  -I_{m} & 0 \\
  0 & I_{m}
\end{array}\right]\end{equation} and \begin{equation}
R^{2}(x)=\left[\begin{array}{cc}
  0 & {\phi}(x) \\
  -{\phi}^{\star}(x) & 2I_{m}
\end{array}\right], \quad 0{\leq }x{\leq}r,\end{equation} where
${\phi}(x){\phi}^{\star}(x)=I_{m}$. From (68) we deduce
that\begin{equation} R(x)=1/2\left[\begin{array}{cc}
  I_{m} &{\phi}(x) \\
  -{\phi}^{\star}(x) & 3I_{m}
 \end{array}\right]\end{equation}
 The matrix $R(x)$ satisfies the following conditions.\\
 1.\emph{The spectrum of $R(x)$ is positive.}\\
 Indeed, we obtain by direct calculation that $[R(x)-I]^{2}=0$.
 Hence the spectrum of the matrix $R(x)$ is concentrated at the point
 ${\lambda}=1$.\\
 2.\emph{The relation \begin{equation}
 jR(x)=R^{\star}(x)j \end{equation} is true.}\\
 It means that $R(x)$ is the j-module of the matrix $W(z)$ which
 satisfies relation (7). From (68) we deduce that \begin{equation}
 R(x)-R^{-1}(x)=\left[\begin{array}{cc}
  -I_{m} & {\phi}(x) \\
  -{\phi}^{\star}(x) & I_{m}
 \end{array}\right].\end{equation}
According to (71) we have
\begin{equation}
 D(x)=j[R(x)-R^{-1}(x)]=\left[\begin{array}{cc}
  I_{m} & -{\phi}(x) \\
  -{\phi}^{\star}(x) & I_{m}
 \end{array}\right].\end{equation} Hence the equality \begin{equation}
 F_{1}(x)=[I_{m}, -{\phi}(x)]
\end{equation} holds.  Using (73) we obtain the
 relations \begin{equation}
 F_{1}(x)jF_{1}^{\star}(x)=0, \end{equation}
\begin{equation}
 F_{1}(x)jF_{1}^{\star}(t)={\phi}(x){\phi}^{\star}(t)-I_{m}
\end{equation}
 Thus in case (69) we deduce from (52) and (74) ,(75) that
 operator the $S_{\xi}$ has the form
 \begin{equation}
 S_{\xi}f=f(x)+\frac{i}{2\pi}P.V.\int_{0}^{\xi}\frac
 {{\phi}(x){\phi}^{\star}(t)-I_{m}}{x-t}f(t)dt .\end{equation}
 The fact that the operator $V$ defined as \begin{equation}
 Vf=\frac{1}{\pi}P.V.\int_{-\infty}^{\infty}\frac{f(t)}{x-t}dt,\quad
 f{\in}L^{2}(-\infty,\infty)\end{equation}
 is unitary implies that \begin{equation}
 S_{\xi}{\geq}0 . \end{equation}\emph{Further we suppose that the
 operator $S_{r}$ is invertible.}\\
 Hence the operators $S_{\xi},\quad \xi{\leq}r$ are invertible as
 well.\\
 \emph{Remark 3}.  If $\phi(x)$ satisfies H\"{o}lder condition
 then there exists such $r>0$ that $S_{r}$ is invertible.\\
Using relation (53) we have
\begin{equation}
 \Phi(x,\xi)+\frac{i}{2\pi}P.V.\int_{0}^{\xi}\frac
 {\phi(x)\phi^{\star}(t)-I_{m}}{x-t}{\Phi}(t,\xi)dt=F_{1}(x) .\end{equation}
 where \begin{equation}
\Phi(x,\xi)=[{\Phi}_{1}(x,\xi),{\Phi}_{2}(x,\xi)].\end{equation}
Here ${\Phi}_{k}(x,\xi)$ are $m{\times}m$ matrix functions
$(k=1,2)$. It follows directly from (73) and (79) that
\begin{equation}
 \Phi_{1}(x,\xi)+\frac{i}{2\pi}P.V.\int_{0}^{\xi}\frac
 {\phi(x)\phi^{\star}(t)-I_{m}}{x-t}{\Phi}_{1}(t,\xi)dt=I_{m},\end{equation}
\begin{equation}
 {\Phi}_{2}(x,\xi)+\frac{i}{2\pi}P.V.\int_{0}^{\xi}\frac
 {{\phi}(x){\phi}^{\star}(t)-I_{m}}{x-t}{\Phi}_{2}(t,\xi)dt=-{\phi}(x),\end{equation}
and \begin{equation}
{\Phi}_{1}(x,\xi){\Phi}_{1}^{\star}(x,\xi)={\Phi}_{2}(x,\xi){\Phi}_{2}^{\star}(x,\xi).
\end{equation} Due to (37) and (54) the formulas \begin{equation}
F_{2}(x)=[-\Phi_{1}(x,1),\Phi_{2}(x,1)],\end{equation}
\begin{equation}
B(\xi)=\frac{1}{2\pi}\int_{0}^{\xi}\left[\begin{array}{cc}
  \Phi_{1}(x,\xi) & \Phi_{2}(x,\xi) \\
 -{\phi}^{\star}(x) \Phi_{1}(x,\xi) &-{\phi}^{\star}(x) \Phi_{2}(x,\xi)
\end{array}\right]dx .\end{equation}are true.\\
\textbf{Example 2}.  We separately consider the partial case of
Example 1 , when $m=1$.\\
It follows from (72) and (82) that \begin{equation}
\Phi_{2}(x,\xi)=-\phi(x)\overline{\Phi_{1}(x,\xi)}.\end{equation}
Hence formula (85) takes the form:
\begin{equation}
B(\xi)=\frac{1}{2\pi}\int_{0}^{\xi}\left[\begin{array}{cc}
  \Phi_{1}(x,\xi) & -\overline{\Phi_{1}(x,\xi)}\phi(x) \\
 -\overline{\phi(x)} \Phi_{1}(x,\xi) &\overline{ \Phi_{1}(x,\xi)}
\end{array}\right]dx .\end{equation}
Comparing formulas (57) and (87) we deduce the representation
\begin{equation}
H(x)=B^{\prime}(x)=a(x)\left[\begin{array}{cc}
  1 & e^{i\alpha(x)}\\
  e^{-i\alpha(x)} & 1
\end{array}\right],\end{equation} where $a(x){\geq}0$,\quad
$\alpha(x)=\overline{\alpha(x)}$. Due to (88) the matrix $jH(x)$
is nilpotent, i.e.
\begin{equation} [jH(x)]^{2}=0 . \end{equation}
\textbf{Example 3}. Let us consider the partial case of Example 1,
 when \begin{equation}
 m=1,\quad \phi(x)=e^{2iux},\quad u=\bar{u}.\end{equation}
 Example 3 plays an important role in the theory of the random
 matrices [4],[7],[14]. Now the operator $S_{\xi}$ takes the form
\begin{equation}
S_{\xi}f=f(x)-\frac{1}{\pi}\int_{0}^{\xi}e^{iu(x-t)}\frac{sinu(x-t)}{x-t}f(t)dt.
\end{equation}The operator $S_{\xi}$ defined by formula (52) is
invertible for all $0<{\xi}<\infty$ (see [4], p.167).\\
We denote by $\Psi(x,\xi,u)$ the solution of the equation
\begin{equation}
\Psi(x,\xi,u)-\frac{1}{\pi}\int_{0}^{\xi}\frac{sinu(x-t)}{x-t}\Psi(t,\xi,u)dt=e^{-iux}.
\end{equation}
Then according to relations (81) and (82) we have \begin{equation}
\Phi_{1}(x,\xi,u)=e^{iux}\Psi(x,\xi,u),\quad
\Phi_{2}(x,\xi,u)=-e^{-iux}\overline{\Psi(x,\xi,u)}.\end{equation}
It follows from (87) and (93), that
\begin{equation}
B(\xi,u)=\frac{1}{2\pi}\int_{0}^{\xi}\left[\begin{array}{cc}
  e^{iux}\Psi(x,\xi,u)& -e^{iux}\overline{\Psi(x,\xi,u)} \\
-e^{-iux}\Psi(x,\xi,u)  & e^{-iux}\overline{\Psi(x,\xi,u)}
\end{array}\right]dx.\end{equation}
 \textbf{Example 4}. Let us
consider the case when $m=1$ and
\begin{equation} J=j=\left[\begin{array}{cc}
  -1 & 0 \\
  0 & 1
\end{array}\right]\end{equation} and \begin{equation}
R(x)=\frac{1}{2}\left[\begin{array}{cc}
  2-|\psi(x)|^{2} & -\overline{\psi(x)}^{2} \\
  {\psi(x)}^{2} & 2+|\psi(x)|^{2}
\end{array}\right], \quad 0{\leq }x{\leq}r.\end{equation}
 The matrix $R(x)$ satisfies the following conditions.\\
 1.\emph{The spectrum of $R(x)$ is positive.}\\
 Indeed, we obtain by direct calculation that $[R(x)-I]^{2}=0$.
 Hence the spectrum of the matrix $R(x)$ is concentrated at point
 ${\lambda}=1$.\\
 2.\emph{The relation \begin{equation}
 jR(x)=R^{\star}(x)j \end{equation} is true.}\\
 It means that $R(x)$ is the j-module of the matrix $W(z)$ which
 satisfies relation (7). From (96) we deduce that \begin{equation}
 R(x)-R^{-1}(x)=jD(x)=jF_{1}^{\star}(x)F_{1}(x),
 \end{equation} where \begin{equation}
 F_{1}(x)=[\psi(x),\overline{\psi(x)}] . \end{equation}
 Using (99) we obtain the
 relations \begin{equation}
 F_{1}(x)jF_{1}^{\star}(x)=0, \end{equation}
\begin{equation}
 F_{1}(x)jF_{1}^{\star}(t)={\psi}^{\star}(x){\psi}(t)-{\psi}(x){\psi}^{\star}(t)
\end{equation}
 Thus  we deduce from (99) and (101) , that the
 operator $S_{\xi}$ in  case (96) has the form
 \begin{equation}
 S_{\xi}f=f(x)+\frac{i}{2\pi}P.V.{\int}_{0}^{\xi}\frac
 {{\psi}^{\star}(x){\psi}(t)-{\psi}(x){\psi}^{\star}(t)}{x-t}f(t)dt .\end{equation}
 \emph{Further we suppose that the operators $S_{\xi}$ are positive and invertible
 $(0<{\xi}{\leq}r)$}.\\
 \emph{Remark 4.}  If $\psi(x)$ satisfies the H\"{o}lder condition,
 then there exists such $r>0$ that the operators $S_{\xi}$ are
 positive and invertible $(0<{\xi}{\leq}r)$.\\
 It follows directly from (99) and (102), that
\begin{equation}
 \Phi_{1}(x,\xi)+\frac{i}{2\pi}P.V.\int_{0}^{\xi}\frac
 {{\psi}^{\star}(x){\psi}(t)-{\psi}(x){\psi}^{\star}(t)}{x-t}{\Phi}_{1}(t,\xi)dt=\psi(x),\end{equation}
\begin{equation}
 {\Phi}_{2}(x,\xi)+\frac{i}{2\pi}P.V.\int_{0}^{\xi}\frac
 {{\psi}^{\star}(x){\psi}(t)-{\psi}(x){\psi}^{\star}(t)}{x-t}{\Phi}_{2}(t,\xi)dt=\overline{{\psi}(x)},\end{equation}
 where \begin{equation}
{\Phi}_{1}(x,\xi){\Phi}_{1}^{\star}(x,\xi)={\Phi}_{2}(x,\xi){\Phi}_{2}^{\star}(x,\xi).
\end{equation} Due to (103) and (104) we have \begin{equation}
F_{2}(x)=[-\Phi_{1}(x,1),\Phi_{2}(x,1)],\quad
\Phi_{1}(x,\xi)=\overline{\Phi_{2}(x,\xi)}
\end{equation}
\begin{equation}
B(\xi)=\frac{1}{2\pi}\int_{0}^{\xi}\left[\begin{array}{cc}
 \overline{ \psi(x)}\Phi_{1}(x,\xi) & \overline{ \psi(x)}\overline{\Phi_{1}(x,\xi)} \\
 {\psi}(x) \Phi_{1}(x,\xi) &{\psi}(x)\overline{\Phi_{1}(x,\xi)}
\end{array}\right]dx .\end{equation}
It follows from (107) that relation (89) is true in this case
as well.\\
\emph{Remark }5. Comparing formulas (69) and (96) we see that
Examples 1 and 3 coincide when $m=1$ and
\begin{equation}
\phi(x)=-\overline {\psi(x)}^{2},\quad |\phi(x)|=1.
\end{equation}
\emph{Remark }6.  If  \begin{equation}
\psi(x)=i\sqrt{\gamma}e^{-iux},\quad 0<\gamma{\leq}1,
\end{equation} due to (96) we have\begin{equation}
R^{2}(x)=\left[\begin{array}{cc}
  1-\gamma & {\gamma}e^{2iux} \\
 {-\gamma}e^{-2iux} & 2+\gamma
\end{array}\right].\end{equation}
The corresponding Riemann-Hilbert problem was considered in
[4].\\
Let us represent $\psi(x)$ in the form  \begin {equation}
\psi(x)=A(x)+iB(x) , \end{equation} where
\begin{equation}
A(x)=\overline {A(x)},\quad B(x)=\overline {B(x)}.\end{equation}
Then the operator $S_{\xi}$ takes the form
\begin{equation}
 S_{\xi}f=f(x)-\frac{1}{\pi}P.V.\int_{0}^{\xi}\frac
 {A(x)B(t)-B(x)A(t)}{x-t}f(t)dt .\end{equation}
 The following partial cases of $\psi(x)$ play an essential role
 in a number of applications [7]:
 \begin{equation}
 \psi_{1}(x)=\sqrt{\pi}[Ai(x)+i{Ai}{\prime}(x)] ,\end{equation}
 where $Ai(x)$ is the Airy function, and \begin{equation}
\psi_{2}(x)=\sqrt{\frac{\pi}{2}}[J_{\alpha}(\sqrt{x})+i\sqrt{x}J_{\alpha}^{\prime}(\sqrt{x})],
\end{equation} where $J_{\alpha}(z)$ is the Bessel function.
\begin{center}{Reference}\end{center}
1. Brodskii M.S., Triangular and Jordan Representation of Linear
Operatorts, \emph{Amer.Math. Soc.1971.}\\
2. Brodskii M.S. and Livsic M.S., Spectral Analysis of
Non-self-adjoint Operators and Intermediate Systems,\emph{ Amer.
Math.
Soc. Transl.(2) 13 , 265-346, 1960}\\
3. Davidson K.R., Nest Algebras,\emph{ Pitnam, Res. Notes Math.,1988}\\
4. Deift P., Its A. and Zhou X, A Riemann-Hilbert Approach to
Asymptotic Problems Arising in the Theory of Random Matrix Models,
and also in the Theory of Integrable Statistical Mechanics,\emph{
Annals
of Math., 146, 149-235, (1997)}.\\
5. Deift P., Integrable Operators,\emph{ Amer. Math. Soc. Transl.
2, v.
189, 69-84, 1999.}\\
6. Gakhov F.D., Boundary Problems,\emph{ Nauka, Moscow, 1977.}\\
7. Harnad J.,Tracy C.A. and Widom H., Hamiltonian Structure of
Equations Appearing in Random Matrices,\emph{ arXiv, 1-18, 1993.}\\
8. Its A.R.,Izergin V.E.,Korepin V.E. and Slavnov N.A., The
quantum Correlation Function as the $\tau$ Function of Classical
Differential Equations, \emph{407-417, Important developments in
soliton theory,
A.S.Fokas and V.E.Zakharov (eds)Springer Verlag , 1993.}\\
9. Larson D. R., Nest Algebras and Similarity
Transformation, \emph{Ann.Math., 125, p.409-427, 1985.}\\
10. Livsic M.S.,Operators, Oscillations, Waves, Open Systems,
\emph{Transl.
of Math. Monographs, 34, Providence, 1973.}\\
11. Potapov V.P.,The Multiplicative Structure of J-contractive
Matrix
Function, \emph{American Math. Society Translation, 15, 131-243, 1960.}\\
12. Sakhnovich L.A., Operators Similar to the Unitary Operator
with Absolutely Continuous Spectrum, \emph{Functional Anal. and
Appl., 2:1,
48-60, 1968}.\\
13. Sakhnovich L.A., Spectral Theory of Canonical Differential
Systems.Method of Operator Identities.\emph{ Operator Theory, Adv.
and
Appl.107, Birkh\"{a}user, 1999.}\\
14. Widom H., Asymptotic for the Fredholm Determinant  of the Sine
Kernel on a Union of Intervals, \emph{ Comm. Math. Phys. 171,
159-180,
1995.}\\
15. Wiener N.,Extrapolation, Interpolation and Smoothing of
Stationary Time Series,\emph{ Cambridge, 1949.}\\
\end{document}